\documentclass[12pt]{article}

\usepackage{amssymb}
\usepackage{amsfonts}

\usepackage{graphicx}
\usepackage{color}
\usepackage{ifpdf}

\usepackage{mathrsfs}
\usepackage{amscd}
\usepackage{amsmath,amsfonts,amssymb,amscd}
\usepackage{indentfirst,graphics,epsfig,psfrag}
\input{epsf}

\usepackage{amsmath,amssymb,amsthm, amscd, amsfonts, graphicx,ifpdf}
\usepackage{lineno}

\newtheorem{thm}{Theorem}[section]
\newtheorem{prop}[thm]{Proposition}

\def\qed{\nopagebreak\hfill{\rule{4pt}{7pt}}
\medbreak}

\def\pf{\noindent {\it Proof.} }

\begin{document}

\title{\Large\bf Note on two results on the rainbow\\ connection number of
graphs\footnote{Supported by NSFC No.11071130.}}
\author{\small Wei Li, Xueliang Li\\
\small Center for Combinatorics and LPMC-TJKLC\\
\small Nankai University, Tianjin 300071, China\\
\small lxl@nankai.edu.cn}
\date{}
\maketitle
\begin{abstract}
An edge-colored graph $G$, where adjacent edges may be colored the
same, is rainbow connected if any two vertices of $G$ are connected
by a path whose edges have distinct colors. The rainbow connection
number $rc(G)$ of a connected graph $G$ is the smallest number of
colors that are needed in order to make $G$ rainbow connected. Caro
et al. showed an upper bound $rc(G)\leq n-\delta$ for a connected
graph $G$ of order $n$ with minimum degree $\delta$ in ``On rainbow
connection, Electron. J. Combin. 15(2008), R57". Recently,
Shiermeyer gave it a generalization that $rc(G)\leq n-
\frac{\sigma_2} 2$ in ``Bounds for the rainbow connection number of
graphs, Discuss. Math Graph Theory 31(2011), 387--395", where
$\sigma_2$ is the minimum degree-sum. The proofs of both results are
almost the same, both fix the minimum degree $\delta$ and then use
induction on $n$. This short note points out that this proof
technique does not work rigorously. Fortunately, Caro et al's result
is still true but under our improved proof. However, we do not know
if Shiermeyer's result still hold.

\end{abstract}

We follow the notation and terminology in \cite{Bondy}, as well as
\cite{Caro} and \cite{Schiermeyer}. The parameter $degree$-$sum$
$\sigma_2(G)$, or simply $\sigma_2$, is defined as $min\{d(u)+d(v) |
u,v\in V(G) \ and \ uv\notin E(G)\}$.

\section{The problem in the proof of Proposition 2.5 of \cite{Caro}}

In \cite{Caro}, Caro et al. got an upper bound of the rainbow
connection number in terms of the minimum degree $\delta$ of a
connected graph, which is stated as follows:
\begin{prop}[Proposition 2.5 of \cite{Caro}]
If $G$ is a connected graph with minimum degree $\delta$, then
$rc(G)\leq n-\delta$.
\end{prop}

In the proof they claimed that they will fix $\delta$ and prove the
proposition by induction on $n$. But the fact is that $\delta_{i}$
of $G_{i}$ is always different from $\delta$. So we cannot use the
induction hypothesis. Fortunately, the proposition is true. We give
an improved proof of the proposition by induction on $n-\delta$.

\pf We prove the proposition by induction on $n-\delta$. The base
case $n-\delta=1$ is trivial since cliques have rainbow connection
number $1=n-\delta$. Assume that the proposition is true for all
connected graphs with $n-\delta < s$. Let us consider a connected
graph with $n-\delta=s$.

Let $K$ be a maximal clique of $G$ consisting only of vertices whose degree is $\delta$.
Since there is at least one vertex with degree $\delta$ and since $G$ is connected
we have $1\leq k=|K|\leq\delta$.

Consider the graph $G'$ obtained from $G$ by deleting the vertices
of $K$. Suppose the connected components of $G'$ are $G_{1}, \ldots,
G_{t}$ where $G_{i}$ has $n_{i}$ vertices and minimum degree
$\delta_{i}$ for $i=1, \ldots, t$. Let $K_{i}\subseteq  K$ be the
vertices of $K$ with a neighbor in $G_{i}$, and assume that
$|K_{1}|\geq |K_{i}|$ for $i=1, \ldots, t$ (notice that it may be
that $t=1$ and $G'$ is connected). Consider first the case that
$K_{1}=K$. We know that $n-n_{i}\geq |K|=k$. For any vertex $v$ in
$G_{i}$, if $d_{G}(v)> \delta$, then $d_{G_{i}}(v)\geq d_{G}(v)-k>
\delta-k$; if $d_{G}(v)= \delta$, then $v$ can not be adjacent to
every vertex in $K$ by the maximality of $K$, so $d_{G_{i}}(v)\geq
\delta-(k-1)>\delta-k$. Thus we have $n-n_{i}\geq
k>\delta-\delta_{i}$, namely $n_{i}-\delta_{i}<n-\delta=s$. By the
induction hypothesis, $rc(G_{i})\leq n_{i}-\delta_{i}$. Clearly, we
may give the edges of $K$ and the edges from $K$ to $G_{1}$ the same
color. Hence,
\begin{align*}
rc(G)&\leq t+\sum_{i=1}^{t}(n_{i}-\delta_{i})=t+n-k-\sum_{i=1}^{t}\delta_{i}\\
&\leq t+n-k-t(\delta-k+1)\\
&=t+n-k-t\delta+tk-t\\
&=n-\delta-(t-1)(\delta-k)\leq n-\delta.
\end{align*}

Now assume that $K_{1}\neq K$ but that $|K_{1}|=k_{1}>1$. Similarly, we have
$n_{i}-\delta_{i}<n-\delta=s$. If there exists a $G_{i}$ with
$\delta_{i}\geq \delta-k+2$, then we can give the edges of $K$ a new color. Hence,
\begin{align*}
rc(G)&\leq t+1+\sum_{i=1}^{t}(n_{i}-\delta_{i})=t+1+n-k-\sum_{i=1}^{t}\delta_{i}\\
&\leq t+1+n-k-(t-1)(\delta-k+1)-(\delta-k+2)\\
&=t+1+n-k-(t-1)(\delta-k)-t+1-\delta+k-2\\
&=n-\delta-(t-1)(\delta-k)\leq n-\delta.
\end{align*}
Otherwise, every $G_{i}$ has $\delta_{i}=\delta-k+1$, namely there
is a vertex $v_{i}$ in $G_{i}$ adjacent to at least $k-1$ vertices
in $K$ for $i=1, \ldots, t$. Since $K_{1}\neq K$, there is a vertex
$u\in K$ but $u \notin K_{1}$. For any $v_{i}$, $i=1, \ldots, t$,
$v_{i}$ is not adjacent to at most $1$ vertex in $K$.
$|K_{1}|=k_{1}>1$, so $v_{i}$ is adjacent to at least $1$ vertex in
$K_{1}$. Now consider any two vertices $v_{i}$ and $v_{j}$, $1\leq i
< j \leq t$. If $v_{i}$ and $v_{j}$ are both adjacent to $u$, then
they have a common neighbor in $K$; if one of them, say $v_{i}$ is
not adjacent to $u$, then $v_{i}$ is adjacent to all vertices in
$K_{1}$, so they also have a common neighbor in $K$. Hence we may
give the edges of $K$ a used color, and so $rc(G)\leq
t+\sum_{i=1}^{t}(n_{i}-\delta_{i})\leq n-\delta$.

Finally, if $k_{1}=1$ (and since $K_{1}\neq K$ we have $k\geq 2$),
contract all of $K$ into a single vertex $v$. Notice that the
contracted graph $G^{*}$ has $n-k+1$ vertices and $\delta^{*}\geq
\delta$. So $n-k+1-\delta^{*}< n-\delta=s$. By induction hypothesis,
$rc(G^{*})\leq n-k+1-\delta^{*}\leq n-k+1-\delta$. Going back to $G$
and coloring the edges of the clique $K$ with another new color, we
obtain $rc(G)\leq n-k-\delta+2\leq n-\delta$. The proof is now
complete. \qed

\section{The problem in the proof of Theorem 9 of \cite{Schiermeyer}}

As a generalization of the above proposition, Schiermeyer got the
following upper bound of the rainbow connection number:

\begin{prop}[Theorem 9 of \cite{Schiermeyer}]
If $G$ is a connected graph with minimum degree-degree $\sigma_2$,
then $rc(G)\leq n- \frac {\sigma_2} 2$.
\end{prop}

In the proof of the result, Schiermeyer used the same proof method
as in Proposition 2.5 of \cite{Caro}. So the proof is also not
correct. Besides, the author claimed that ``each pair of nonadjacent
vertices of $G_{i}$ has degree-sum at least $\sigma_{2}(G)-2(k-1)$
in $G_{i}$'', this is not correct. We can give a counter-example.
Just take $t$ copies of $K_{\delta+4}$ as $H_{1}, \ldots, H_{t}$,
$t< \frac{\delta+1}{2}$. $v_{i,1}$ and $v_{i,2}$ are two nonadjacent
vertices, $i=1, \ldots, t$. Connect $v_{i,1}$ and $v_{i,2}$ to any
$2t$ vertices in $H_{i}$, respectively, the neighbor of $v_{i,1}$
and $v_{i,2}$ can be same or different. Let $K=K_{\delta-2t+1}$ and
connect $v_{i,1}$, $v_{i,2}$ to all vertices in $K$. Denote the
graph by $G$. Then the vertices in $H_{i}$ have degree at least
$\delta+3$, $v_{i,1}$, $v_{i,2}$ have degree $\delta+1$ and the
vertices in $K$ have degree $\delta$. Since every vertex with degree
$\delta$ is adjacent to every vertex with degree at most $\delta+1$
and other vertices have degree at least $\delta+3$,
$\sigma_{2}(G)=d(v_{i,1})+d(v_{i,2})=2(\delta+1)$. In
$G_{i}=G[H_{i}\cup\{v_{i,1}, v_{i,2}\}]$, $v_{i,1}$ and $v_{i,2}$
are two nonadjacent vertices and
$d_{G_{i}}(v_{i,1})+d_{G_{i}}(v_{i,2})=4t$. But
$\sigma_{2}(G)-2(k-1)=2(\delta+1)-2(\delta-2t+1-1)=4t+2>4t=d_{G_{i}}(v_{i,1})+d_{G_{i}}(v_{i,2})$,
contradicting the claim. On the other hand, since $d_{G_{i}}(v)\geq
d_{G}(v)-k$, we can get $s_{i}\geq \sigma_{2}(G)-2k$. Therefore
\begin{align*}
rc(G)&\leq t+\sum_{i=1}^{t}(n_{i}-\frac{s_{i}}{2})=t+n-k-\sum_{i=1}^{t}\frac{s_{i}}{2}\\
&\leq t+n-k-t(\frac{\sigma_{2}(G)}{2}-k)\\
&=t+n-k-t\frac{\sigma_{2}(G)}{2}+tk\\
&=n-\frac{\sigma_{2}(G)}{2}-(t-1)(\frac{\sigma_{2}(G)}{2}-k)+t\\
&\leq n-\frac{\sigma_{2}(G)}{2}+t.
\end{align*}
When $t$ is large, this bound is far from
$n-\frac{\sigma_{2}(G)}{2}$ given by the author of
\cite{Schiermeyer}. However, we do not know whether his theorem is
still true.

\end{document}